\documentclass{amsart}
\input{psfig.sty}

\usepackage{amssymb}
\usepackage{hyperref}

\def\abs#1{\left|#1\right|}

\newtheorem*{conjecture}{Conjecture}

\begin{document}

\title{An algorithm to approximate reals by rationals of the form $a/b^2$. }
\author{I. Jim\'enez Calvo}


\address{Grupo GASS-An\'alisis, Seguridad y Sistemas. Facultad de Inform\'atica, U.C.M. Ciudad Universitaria, 28040-Madrid, Spain.}
\email{ijcalvo@terra.es}
\date{14th february 2007}

\maketitle

\begin{abstract}
This paper develops an algorithm to search for approximations to a real $\xi$ which finds rationals of the form $a/b^2$, such that $\abs { \xi - \frac a {b^2} } < \frac {c(\xi)} { b^3}$, in $B^{1/2}$ polynomial time operations, for $b \le B$. An heuristic approach and computational data give support to the conjecture of O. Karpenkov stating that there are infinitely many solutions in integers $a$ and $b$ for such approximation.
\end{abstract}

\section{Introduction}
Recently, O. Karpenkov \cite{OK} posed the problem of the approximation of reals by rationals of the form $a/b^2$ as a natural subject of number theory and for its incidence in the solution of cases of the Shr\"odinger equation. He presents a proved lower estimate for the approximation of a real $\xi$ as 
$$
\abs { \xi - \frac a {b^2} } < \frac c { b^3\ln ^{1+\varepsilon} b},
$$
and conjectures that there exists a constant $c(\xi)$ such that the inequality
$$
\abs { \xi - \frac a {b^2} } < \frac {c(\xi)} { b^3},
$$
has infinitely many solutions in integers $a$ and $b$.

In this paper, we study the problem from an heuristic point of view and develop an algorithm that can compute an approximation of the form $a/b^2$ in $\sqrt b$ polynomial time operations. The heuristic approximation to the problem and the computational data obtained support the conjecture of O. Karpenkov.

\section{An heuristic approach}
Following a brute force scheme to search for approximations of the form $a/b^2$ to an irrational $\xi$, we can compute, for $b=1,2,\cdots$, the optimal value of $a$ as the nearest integer to $b^2\xi$. For each $b$ value we get an approximation
$$
\abs{\xi - \frac a {b^2}} = \frac \delta {b^2},\quad 0< \delta < \frac 1 2.
$$
It seems reasonable to suppose that the values of $\delta$ are distributed uniformly in this interval. The probability of getting a value $\delta < c/b^{1+\varepsilon}$, where $c$ is a constant, will be, under that supposition, $p(b)= 2c/b^{1+\varepsilon}$. Note that, in this case,
$$
\abs{\xi - \frac a {b^2}} < \frac c {b^{3+\varepsilon}}.
$$
Moreover, we can estimate the number of rationals with the above approximation as
$$
\sum_{b=1}^\infty p(b)=2c\sum_{b=1}^\infty \frac 1 {b^{1+\varepsilon}}.
$$ 
Since the summatory converges for $\varepsilon>0$, we can guess that there may be only finitely many rational approximations of the form $a/b^2$ of any real $\xi$, for $\varepsilon>0$.

When the exponent in $b$ is exactly 3, the summatory is the harmonic series and
$$
\sum_{b=1}^B \frac 1 b = \ln B + \gamma + O(\frac 1 B),
$$
where $\gamma= 0.57721\cdots$ is the Euler's constant (see \cite{LeVeque} Theorem 6.10, p. 136). Therefore, the number of approximations with $b<B$ may be estimated as $2c(\gamma+ \ln B)$. Note that the fractions of that form are not necessarily irreducible in this estimation.

\section{The algorithm}
We consider the following modular quadratic equation,
\begin{equation}\label{eqmod}
Pb^2 \equiv \alpha \pmod Q.
\end{equation}
Since $Pb^2=aQ+\alpha$, for some integer $a$, we can write
$$
\frac P Q = \frac a {b^2} + \frac \alpha {b^2Q}.
$$
Let $\xi$ be a real number. From the equation
$$
\xi - \frac a {b^2} = \xi -\frac P Q + \frac \alpha {b^2Q},
$$
we have that
\begin{equation}\label{bound}
\abs{\xi - \frac a {b^2}} \le \abs{ \xi -\frac P Q} + \abs{\frac \alpha {b^2Q}}.
\end{equation}
If we choose $P/Q$ as a convergent of $\xi$ computed from its continued fraction expansion, we have that 
\begin{equation}\label{aprox}
\abs{\xi-\frac P Q} < \frac 1 {\sqrt{5}Q^2}
\end{equation}
from the well known result due to Hurwitz (see, by example Theorem 193, p. 164 in \cite{H&W}).

Suppose that we are able to find a solution to equation (\ref{eqmod}) such that $b \le kQ^{2/3}$ and $\abs \alpha \le tQ^{1/3}$ for certain constants $k$ and $t$. Then we have that
\begin{equation}\label{error}
\abs{\xi - \frac a {b^2}} <  \left(  \frac {k^3} {\sqrt 5} + kt   \right) \frac 1 {b^3}.
\end{equation}

If we solve equation (\ref{eqmod}) for some value of $\alpha$, we have a probability $kQ^{-1/3}$ of getting a value $b\le Q^{2/3}$, under the supposition that the quadratic residues are distributed at random. Since we check $2Q^{1/3}$ values of $\alpha$, we have the chance of getting $2tk$ values of $b$ in average which gives approximations of $\xi$ that fulfill (\ref{error}). For this task, we must solve $2tQ^{1/3}$ quadratic equations for which there are polynomial time algorithms whenever the factorization of $Q$ is known. Suppose that a value of $b$ is obtained such $b=Q^{2/3}$, so the number of quadratic equations solved to obtain it is $O(b^{1/2})$. Take into account that there are methods of factorization that can factorize $Q$ in $Q^u$ steps with $u$ well below $1/3$ (see \cite{Riesel} p. 218), then the complexity of this stage is smaller than the solution of the quadratic equations.

Using the described method, there may be approximations fulfilling (\ref{error}) that scape the search of the algorithm. First, a convergent $P/Q$ may have an approximation to $\xi$ better than (\ref{aprox}), mainly for small $Q$, what allows that higher values of $b$ keep bounded the first summand in (\ref{bound}). Besides, both summands can compensate if they are of opposite sign. For these reasons, in practice, it is better to use a brute force scheme for the fist values of $b$ (say for $b<1000$) and then to apply the described algorithm. Likewise, it is better to check the values of $\abs \alpha \le Q^{1/3+\varepsilon}$ (by example of $\abs \alpha \le Q^{0.35}$). Since, once computed the roots of (\ref{eqmod}), checking the values of $b$ is less time consuming, it is better to rise the bound of $b$ to $b\le Q^{3/4}$. With these trimmings, the algorithm finds all items existing (checked by brute force) at least for $b < 2*10^9$.

The algorithm was programmed using the PARI package \cite{PARI} and compiled with the utility \textit{gp2c}. The executable was launched in an AMD64 +3200 computer over one week, searching for rational approximations of $\pi$ such that the difference is less than $1/b^3$ in absolute value.  There were found 94 rational approximations of the prescribed form which can be found in \cite{pi} (note that only $b$ values are tabulated because the numerator $a$ can be computed from it). The higher rational approximation found correspond to 
$$
a=36266840658555398816245943123914613560.
$$
$$
b=3397660065732068041.
$$
The following figure plots the number of items versus the expected value $1+2\ln b$ showing that both agree quite well. Another irrationals as $e$, $\sqrt 2$, the golden mean $(1+\sqrt 5 )/2$, even numbers suspicious of irrationality as the Euler's constant were tested with the same results. 
\begin{figure}
\centerline{\psfig{figure=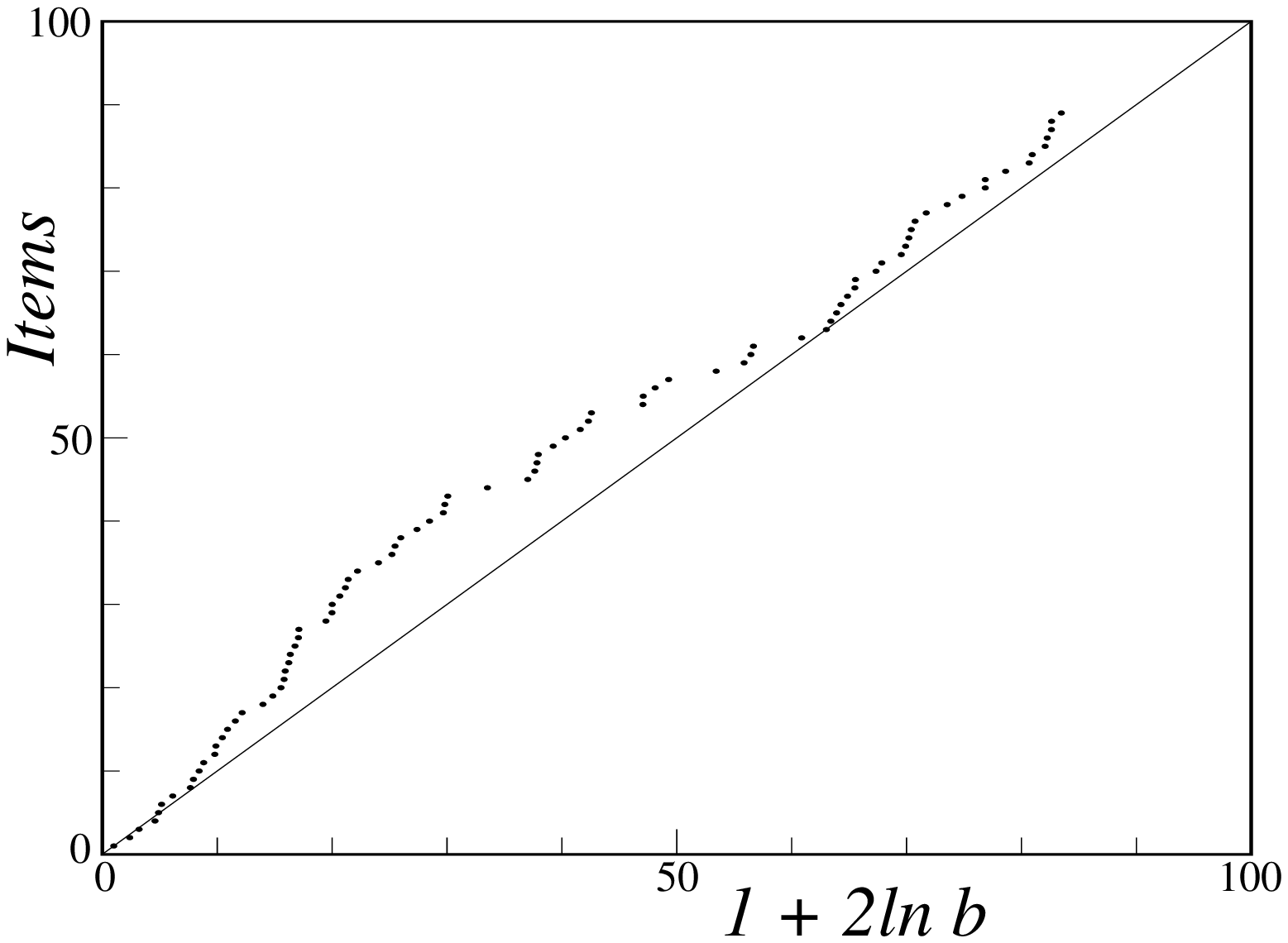,height=8cm}}
\end{figure}

\section{Approximating irrationals by rationals of the form $a/p$, $p$ prime.}
O. Karpenkov, in his cited paper, also raised the problem of approximating a real with a rational with prime denominator. We follow the procedure described in the preceding section, to find approximations of that kind. Consider the linear modular equation 
\begin{equation}\label{eqmod2}
Px \equiv \alpha \pmod Q.
\end{equation}
Since $Px= aQ+\alpha$
\begin{equation}\label{bound2}
\abs{\xi - \frac a x} \le \abs{ \xi -\frac P Q} + \abs{\frac \alpha {xQ}}.
\end{equation}
We take $P$ and $Q$ from a convergent of $\xi$ and solve the equation (\ref{eqmod2}) for the lowest absolute values of $\alpha$. Since the density of primes is about $1/\ln n$, we may expect that we will find a $x$ prime for $\alpha$ of the order of $ln Q$ and let that prime be $p$. This constructive approach makes reasonable the
\begin{conjecture}
For any real $\xi$, there exists a constant $c(\xi)$ such that the inequality
$$
\abs{\xi - \frac a p} < \frac {c(\xi)\ln p} {p^2},
$$
with $p$ prime, has infinitely many solutions.
\end{conjecture}

\vskip 3mm
{\bf Dedicatory:} I dedicate this paper to the memory of my friend Manuel Torres Hernanz. Since my first collaboration with him, two decades ago, computing the Fourier transform of his magnetic measures till a recent talking about the different ways to get rational approximations to $\sqrt 2$, I have enjoyed his knowledge, creativity, and humanistic spirit. I miss him and I am very indebted to his generosity.

\end{document}